Galina I. Sinkevich

https://orcid.org/0000-0002-8595-1686

Saint Petersburg State University

of Architecture and Civil Engineering,

(Saint-Petersburg, Russia)

### On the history of negative and complex numbers' interpretation

*Abstract*. The history of the development of the concept of complex numbers from the 16th to 19th centuries. The origin and refinement of the geometric and physical meaning of complex numbers, the emergence of vectoral analysis.

The development of the concept of numbers has a long history. In antiquity, only natural numbers were called numbers: 1, 2, 3, …, and aliquot fractions and proportions (ratios of natural numbers) were known. In Ancient Greece, it was discovered that values existed which could not be expressed in any ratio, for example the diagonal of the square of one. They were calculated approximately, using the method of exhaustion. The coefficients of equations, the roots of the equations could only be positive. If in solving tasks, a negative value arose, it was considered to be meaningless, as this quantity could not be less than nothing. Sometimes people said that in the commercial sense it could mean a debt. In Russian in the 18th century negative numbers were called losses[1]. So in solving equations, only positive roots were looked for. This was how things were until the Renaissance.

A negative number was also not recognized as a proper mathematical entity because for values, the rule of ratio had to be carried out: if the left part is a ratio that is less than the greater one, then the right part of the ratio should also be a ratio of the lesser to the greater. But for ratios of the kind $\frac{-1}{1} = \frac{1}{-1}$ this was not carried out. Additionally, a quantity could not be less than nothing, whether this related to natural or rational numbers.

In 1494, Luca Pacioli (1445–1517) wrote the treatise *Summa de arithmetica, geometria, proportioni et proportionalità* (The sum of arithmetic, geometry, proportions and proportionality), in which he collected the knowledge of arithmetic held by Europeans and Indians.

---

[1] Euler, L. Universal'naya arifmetika g. Leongarda Euler'a. Perevedennaya s nemeckogo podlinnika studentami Petrom Inohodcovym i Ivanom Yudinym. Tom 1, soderzhashchij v sebe vse obrazy algebraicheskogo vychisleniya. Saint Petersburg: Imp. academiae scientiarum Petropolitanae, 1768. - p. 10.

In 1544, Michael Stifel was the first[2] to state that negative numbers were numbers less than zero (below zero)[3]. Since this time, the concept of the numerical scale gradually began to form, in which positive numbers are located to the right of zero in increasing order, and negative numbers to the left of zero in decreasing order. In the 17[th] century, this was reflected in the chronological scale – for the first time, various chronologies were given in scale, which had a point of reference, a positive direction (AD) and a reverse calculation (BC) – Joseph Juste Scaliger (1540–1609) and Dionysius Petavius (1583–1652). In 1742, the Swedish astronomer Anders Celsius (1701—1744) created a temperature scale with zero as the point of reference. The scale received its modern form thanks to Carl Linneus[4].

Before the 17[th] century, movement was only examined as even directly or in a circle. More complex trajectories of movement emerged in the 17[th] century. It became possible to describe uneven movement after the creation of mathematical analysis. Mechanical interpretation of a negative number was first devised by John Wallis, who described the example of forward movement in a straight line, 5 yards forwards, and then 8 yards backwards[5].

## 1545  Gerolamo Cardano's *Artis Magnæ, Sive de Regulis Algebraicis Liber Unus* (Ars Magna, The Great Art)

1545 is considered to be the year in which complex numbers were discovered. Gerolamo Cardano, studying the solution of a cubic equation, in intermediary insertions came across the case of false roots of an auxiliary equation, which were subsequently eliminated. He only searched for positive roots, and called negative roots impossible, and roots from negative values truly sophistic[6]. Remember that in these times, there was no algebraic symbols, or formulas. Rules were expressed in words. Here is page 287 from Cardano's *Ars Magna* and a translation of it: Figure 1. Page 287 from Cardano's *Artis Magna*

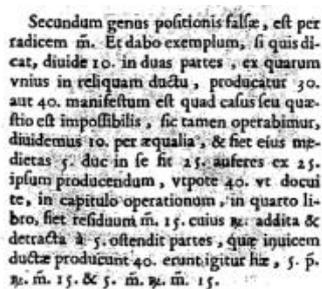

Cardano, *Ars Magna*, Chapter XXXVII (De regula falsum ponendi – the rule of a false proposition, a negative integer): "The second type of false solution lies in the root of a negative quantity (per radicem $m$). I will give an example. If someone demands to divide 10 into two parts, which when multiplied would equal 30 or 40, it is obvious that this case or question is impossible. But we will do as follows: we divide 10 in half, half will be 5; multiplied by itself, it comes to 25. Then deduct from 25 what should result from multiplication, say 40 – as I explained to you in the chapter on operations in the $4^{\text{th}}$ book; then we are left with $m:15$; if we deduct from this $R$ and add to 5 or deduct from 5, then we get parts which when multiplied together give 40. Thus, these parts will be $5p:Rm:15$ and $5m:Rm:15$ ".

Here Cardano examines the equation $x(10-x)=40$, or $x^2+40=10x$ (the equation here should be written so that the coefficients are positive). Cardano solves this according to his rule, in modern notation: $x=5\pm\sqrt{25-40}=5\pm\sqrt{-15}$. The colon was used instead of a full stop. Cardano shows that the derivation of the roots is 40. In the text $p$ means plus, and $m$ means minus. $R$ means radix, root; $5p:Rm:15$ means $5+\sqrt{-15}$; $5m:Rm:15$ signifies $5-\sqrt{-15}$; $25m:m:15$ quod est 40 signifies $25-(-15)=40$. Cardano uses the fact that the complex roots when multiplied give a real number. The history of the discovery of the formula for solving cubic equations is described in the books: Niccolò Tartaglia. *Quesiti et inventioni diverse, dialogo con interlocutori principali Francesco Maria della Rovere e Gabriele Tadino e argomenti diversi: aritmetica, geometria, algebra, statica, topografia, artiglieria, fortificazioni, tattica*. 1546.; Bortolotti, E. *La storia della matematica nella Università di Bologna by Ettore Bortolotti*. Bologna: N. Zanichelli, 1947. 226 p.; Guter R., Polunov Yu. Girolamo Kardano. M.: Znanie, 1980. 192 p. (Гутер Р., Полунов Ю. *Джироламо Кардано*. М.: Знание, 1980 г. 192 с.); S.G. Gindikin. *Rasskazy o fizikah i matematikah* (izdanie tret'e, rasshirennoe). M.: MCNMO, NMU, 2001. 448 p. (С.Г. Гиндикин. *Рассказы о физиках и математиках* (издание третье, расширенное). М.: МЦНМО, НМУ, 2001 г. 448 с.).

### 1572 Rafael Bombelli's *Algebra*

In 1572, a follower of Cardano, the hydraulic engineer Rafael Bombelli (1526–1572), wrote the book *Algebra*[7], where he first introduced the rules of arithmetic operations on negative numbers, and examined the solution of the cubic equations with roots of negative values. In solving these equations, where in auxiliary equations under the sign of the cubic radical it was possible to select the cube of the sum or the difference, and thus extract the cubic root, Bombelli showed that roots of negative values are mutually destroyed, as the components are mutually conjugate. Bombelli showed the possibility of determining the ratio of equality, sum and derivation of complex numbers. But roots of negative values still had no physical or geometric meaning. Bombelli, as a hydraulic engineer, saw them as a useful auxiliary construction.

---

[7] Bombelli R. L'Algebra opera. Divisa in tre libri. Bologna: Nella stamperia do Guovanni Rossi. 1572.

### 1637 René Descartes' *Geometry*

In 1637, René Descartes (1595–1650) published his *Geometry*, in which he called false roots "imaginary" (imaginariae). The term "real root" first appeared there. He called negative roots false roots. Descartes devised imaginary roots in solving the problem of crossing a circle with a parabola. Descartes examined cases of their intersection, touching, and the case "when the circle does not cross the parabola at any point, and this means that the equation does not have true or false roots, and that they are all imaginary". "There is no value which corresponds to these imaginary roots"[8].

### 1685 John Wallis' *Algebra*

The first mathematician to attempt to give a geometric and physical interpretation of negative and imaginary numbers was John Wallis. This was in 1685, in his treatise *Algebra*. He explains negative numbers in a problem of displacement: "Yet is not that Supposition (of Negative Quantities,) either Unuseful or Absurd4 when rightly understood. And though, as to the bare Algebraick Notation, it import a Quantity less than nothing: Yet, when it comes to a Physical Application, it denotes as Real a Quantity as if the Sign were +4 but to be interpreted in a contrary sense.

As for instance: Supposing a man to have advanced or moved forward, (from *A* to *B*) 5 Yards; and then to retreat (from *B* to *C*) 2 Yards: If he asked, how much he had Advanced (upon the whole march) when at *C*? or how many Yards he is now Forwarder than when he was at *A*? I find (by Subducting 2 from 5,) that he is Advanced 3 Yards. (Because +5 −2=+3.)

But if, having Advanced 5 Yards to *B*, he thence Retreat 8 Yards to *D*; and it be asked, How much he is Advanced when at *D*, or how much Forwarder than when he was at *A*: I say −3 Yards. (Because +5 −8= −3.) That is to say, he is advanced 3 Yards left than nothing.

Which in propriety of Speech, cannot be, (since there cannot be less than nothing.) And therefore as to the Line *AB Forward*, the case is Impossible.

But if (contrary to the Supposition,) the Line from *A*, be continued *Backward*, we shall find *D*, 3 Yards behind *A*. (Which was presumed to be *Before* it.)

---

And thus to say, he is Advanced −3 Yards; is but what we should say (in ordinary form of Speech,) he is *Retreated* 3 Yards; or he wants 3 Yards of being so Forward as he was at *A*"[9]. You can see this text in the image below, the spelling and italics are original.

Figure 2. J. Wallis, *Algebra*, page 265.

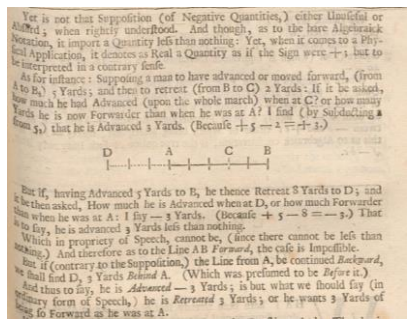

Imaginary numbers are like sides of a lost square field of earth[10]. An imaginary value for him is the "middle proportional value between a positive and negative value". In his drawing we see that the imaginary number is a section of the *BP*. This is Wallis's argument:

"What hath been already said of $\sqrt{-bc}$ in Algebra, (as I Mean Proportional between a Positive and a Negative Quantity:) may be thus Exemplified in Geometry.

If (for instance,) Forward from *A*, I take *AB*=+*b* and Forward from thence, *BC*=+*c*; (making $AC = +AB + BC = +b + c$, the Diameter of a Circle:) Then is the Sine[11], or Mean Proportioned $BP = \sqrt{+bc}$.

But if Backward from *A*, I take $AB = -b$; and then Forward from that *B*, $BC = +c$; (making $AC = -AB + BC = -b + c$; the Diameter of the Circle:) Then is the Tangent or Mean Proportional $BP = \sqrt{-bc}$.

So that where $\sqrt{+bc}$ signifies a Sine; $\sqrt{-bc}$ shall signify a Tangent, to the same Arch (of the same Circle,) *AP*, from the same point *P*, to the same Diameter *AC*.

Suppose now (for further Illustration,) A Triangle standing on the Line *AC* (of indefinite length;) whose one Leg *AP*=20 is given; together with (the Angle *PAB*, and consequently) the Height *PC*=12; and the length of the other Leg *PB*=15: By which we are to find the length of the Base *AB*.

'Tis manifest that the Square of *AP* being 400; and of *PC*, 144; their Difference 256 (=400 −144) is the Square of *AC*.

And therefore AC( $= \sqrt{256}$ )=+16, or −16; Forward or backward according as we please to take the Affirmative or Negative Root. But we will here take the Affirmative.

Then, because the Square of *PB* is 225; and of *PC*, 144; their Difference 81, is the Square of *CB*. And therefore $CB = \sqrt{81}$; which is indifferently, +9 or –9; And may therefore be taken Forward or Backward from *C*. Which gives a Double value for the length of *AB*; to wit, $AB = 16 + 9 = 25$; or $AB = 16 - 9 = 7$. Both Affirmative. (But if we should take, Backward from *A*, *AC*= −16; then *AB*= −16+9= −7, and *AB*= −16 −9= −25. Both Negative.)

Figure 3. J. Wallis, *Algebra*, page 266.

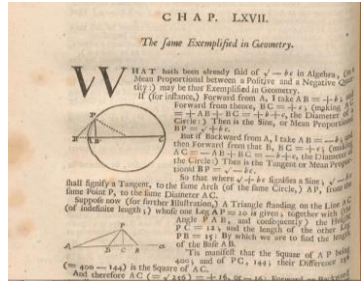

Suppose again, *AP*=15, *PC*=12, (and therefore $AC = \sqrt{:225 - 144:} = \sqrt{81} = 9$,) *PB*=20 (and therefore $BC = \sqrt{:400 - 144:} = \sqrt{256} = +16$ or −16:) Then is *AB*=9+16=25, or $AB = 9 - 16 = -7$. The one Affirmative, the other Negative. (The same values would be, but with contrary Signs, if we take $AC = \sqrt{81} = -9$: That is, $AB = -9 + 16 = +7$, $AB = -9 - 16 = -25$.)

In all which cases, the Point *B* is found, (if not Forward, at least Backward,) in the Line *AC*, as the Question supposeth.

And of this nature, are those Quadratick Equations, whose Roots are Real, (whether Affirmative or Negative, or partly the one, partly the other;) without any other Impossibility than (what is incident also to Lateral Equations,) that the Roots (one or both) may be Negative Quantities.

But if we shall Suppose, *AP*=20, *PB*=12, *PC*=15, (and therefore $AC = \sqrt{175:}$) When we come to Subtract as before, the Square of *PC* (225,) out of the Square *PB* (144,) to find the Square of *BC*, we find that cannot be done without a Negative Remainder, 144 −225= −81.

Figure 4. J. Wallis, *Algebra*, page 267.

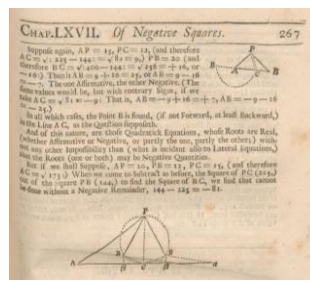

So that the Square of *BC* is (indeed) the Difference of the Squares of *PB*, *PC*; but a defective Deference; (that of *PC* proving the greater, which was supported the Lesser; and the Triangle *PBC*, Rectangled, not as supported at *C*, but at *B*:) And therefore $BC = \sqrt{-81}$

Which gives indeed (as before) a double value of *AB*, $\sqrt{175}, +\sqrt{-81}$, and $\sqrt{175}, -\sqrt{-81}$: But such as requires a new Impossibility in Algebra, (which in Lateral Educations doth not happen;) not that of a Negative Root, or a Quantity less than nothing; (as before,) but the Root of a Negative Square. Which in strictness of speech, cannot be: since that no Real Root (Affirmative or Negative,) being Multiplied into itself, will make a Negative Square.

This Impossibility in *Algebra*, argues an Impossibility of the case proposed in Geometry; and that the Point *B* cannot be had, (as was supported,) in the Line *AC*, however produced (forward or backward,) from *A*.

Yet are there Two Points designed (out of that Line, but) in the same Plain; to either of which, if we draw the Lines *AB*, *BP*, we have a Triangle; whose Sides *AP*, *PB* are such as were required: And the Angle *PAC*, and Altitude *PC*, (above *AC*, though not above *AB*,) such as proposed; And the Difference of Squares of *PB*, *PC*, is that of *CB*.

And like as in the first case, the Two values of *AB* (which are both Affirmative,) make the double of AC, $(16 + 9, + 16 - 9, = 16 + 16 = 32:)$ So here, $\sqrt{175} + \sqrt{-81}, + \sqrt{175} - \sqrt{-81} = 2\sqrt{175}.$

And (in the Figure,) though not the Two Lines themselves, *AB*, *AB*, (as in the First case, where they lay in the Line *AC*;) yet the Ground-lines on which they stand, *A*β, *A*β, are Equal to the Double of *AC*: That is, if to either of those *AB*, we join *B*α, equal to the other of them, and with the same Declivity; *AC*α (the Distance of *A*α) will be a Streight Line equal to the double of *AC*; as is AC α in the First case.

The greatest difference is this; That in the first Case, the Points *B*, *B*, lying in the Line *AC*, the Lines *AB*, *AB*, are the same with their Ground-Lines, but not so in this case, where *BB* are so raised above ββ (the respective Points in their Ground-Lines, over which they stand,) as to make the case feasible; (that is, so much as is the versed Sine of *CB* to the Diameter *PC*:) But in both *AC*α (the Ground-Lines of *AB*α) is Equal to the Double of *AC*.

So that, whereas in case of Negative Roots, we are to say, The Point *B* cannot be found, so as is supported in *AC* Forwards, but Backwards from *A* it may in the same Line: We must here say, case of a Negative Square, the Point *B* cannot be found so as was supported, in the Line *AC*; but Above that Line it may in the same Plain[12]."

Unfortunately, Wallis's suggestion that complex number were not located on a straight lane, but on a complex plane, was incomprehensible to his contemporaries. For a number to become a mathematical object, its ratios had to be determined (equality, less, more, i.e. order) and operations on objects. But it was unclear whether an operation on complex numbers would always lead to a number of this kind, i.e. $x + y\sqrt{-1}$.

[12] Ibid, p. 266–268.

### 1702, Leibniz

Gottfried Wilhelm Leibniz tried to prove this in 1702, but failed. In the article *Graphic proof of the new analysis for recognizing infinity in relation to sums and quadratics* [13], breaking down the binomial into multiples, Leibniz arrived at the result

$$x^4 + a^4 = \left( x + a\sqrt{\sqrt{-1}} \right)\left( x - a\sqrt{\sqrt{-1}} \right)\left( x + a\sqrt{-\sqrt{-1}} \right)\left( x - a\sqrt{-\sqrt{-1}} \right),$$ and concluded that an

imaginary number of a different kind existed. He called imaginary numbers idealis mundi monstro: "Itaque elegans et mirabile effugium repetir in illo Analyseos miraculo, idealis mundi monstro, pene inter Ens et non-Ens Amphibio, quod radicem imaginariam appellamus[14]". – "What we call the imaginary root is an elegant and wonderful invention in this incredible analysis, the prototype of a monster of the world, an amphibian between being and non-being".

### 1712 The logarithm of the negative and imaginary number

Before 1702, imaginary numbers were only seen as roots of negative values. In 1702, Johann Bernoulli encountered the problem of calculating the logarithm of a complex number. By 1712, Bernoulli and Leibniz had argued about what the logarithm of a negative number was. For the positive number *a*, it was fair to state $\ln \sqrt{a} = \frac{1}{2}\ln a$. Continuing the argument, it may be concluded that $\ln i = \ln \sqrt{-1} = \frac{1}{2}\ln(-1)$. But what is equal to $\ln(-1)$? Leibniz believed that it must be complex. Bernoulli, and d'Alembert after him, believed that it was substantial. The English mathematician and astronomer R. Cotes, in his work *Logometria*, 1714, published in Philosophical Transactions in 1717, placed the formula $\ln(\cos x + i\sin x) = xi$, expressed in the following words: "If any arc of a quarter of a circle, described by the radius, *CE*, has the sine *SX* and the cosine to the quarter *XE* and if the radius *CE* is taken for the module, then the arc will be the measure of the ratio $EX + XC\sqrt{-1}$ to *CE*, multiplied by $\sqrt{-1}$". Cotes did not give any applications for this.

In 1749 Euler proved it, confirming that Leibniz was correct. Now we know this formula as $\operatorname{Ln} z = \ln|z| + i\varphi + 2k\pi i$.

---

### 1707 and 1722. Abraham de Moivre's trigonometric representation

In 1707, and later in 1722, Abraham de Moivre made a trigonometric interpretation of the complex number. Cubic equations and higher were solved not only by the algebraic method, but also by the trigonometric method, using the sinus of short arcs[15]. In 1594, Francois Viète had solved an equation of the $45^{th}$ degree using this method. Using known ratios, Moivre reached a formula of raising the degree and extracting the root of a natural degree (up to the $7^{th}$) of a complex number. It is interesting that he examined the circular arc $x^2 + y^2 = 1$, and then the hyperbolic arc $x^2 - y^2 = 1$, which led him to the idea of the imaginary substituting $y = v\sqrt{-1}$ [16]. But even when he presented a complex number in trigonometric form, Moivre did not depict it on a surface.

### L. Euler

In Petersburg in the 1730s–1740s, Leonhard Euler developed the rudiments of the theory of functions of the complex variable. In his works, Euler moved from coordinates of a point $(x, y)$ to the complex number $p = x \pm \sqrt{-1} y$, and represented it in the polar coordinates $p = s\left(\cos \omega \pm \sqrt{-1} \sin \omega\right)$. This concept was used after Euler by Lagrange and other mathematics in two-dimensional problems of mathematical physics, but at that time there was no geometric, let alone physical concept of operations on complex numbers. In 1743, Euler created the method of solving linear differential equations of higher orders, in which in solving characteristic algebraic equations, imaginary numbers arise. At the same time, the general solution of the equation is valid[17].

In 1748, Euler proved Moivre's formula for all valid *n*. Now it is proven as a consequence from Euler's formula $e^{i\varphi} = \cos \varphi + i \sin \varphi$. Euler published this formula in an article of 1740 and in the $7^{th}$ chapter of his book *Introduction to an analysis of infinitesimals* (Introductio in analysin infinitorum, 1748 г.)[18]

Euler gradually gained an understanding of the concept of the complex number. A large number of observations of his own mathematical studies did not always find a geometric or physical interpretation. A.I. Markushevich noted this fact. In 1741 Euler in a letter to Goldbach

(9.XII. 1741) wrote: "I recently found a wonderful paradox, that this value of the expression

$\dfrac{2^{+\sqrt{-1}} + 2^{-\sqrt{-1}}}{2}$ is very close to $\dfrac{10}{13}$, and this fraction only differs from reality by a few millionths.

The true value of this expression is the cosine of the arc $0,6931471805599...$". The meaning of this statement becomes clear if we represent $\dfrac{2^{+\sqrt{-1}} + 2^{-\sqrt{-1}}}{2}$ in the form $\dfrac{e^{+\sqrt{-1}\ln 2} + e^{-\sqrt{-1}\ln 2}}{2}$, which according to Euler's formula from "An introduction to an analysis of infinitesimals",

$\cos v = \dfrac{e^{+v\sqrt{-1}} + e^{-v\sqrt{-1}}}{2}$ gives $\cos \ln 2$. The value of $\ln 2$ is given by Euler[19].

### 1749/51, Euler

In his article *Studies on imaginary roots of equations* (Recherches sur les racines imaginaires des équations. Mém. Ac. Berlin, (1749) 1751) examined the issue of the possible form of a complex number. "A quantity is named imaginary is it is not greater than zero, not less than zero or equal to zero; this, accordingly, is something impossible, for example $\sqrt{-1}$ or even $a + b\sqrt{-1}$, as this quantity is not positive, or negative, or zero"[20]. Euler examines the main theorem of algebra that he has proven as a separate case of the following proposal "any imaginary quantity is always formed by two members, one of which is a real quantity indicated by *M*, and the other is a derivate of this real quantity *N* by $\sqrt{-1}$; thus $\sqrt{-1}$, is the only source of all imaginary expressions"[21].

For proof, Euler applied to numbers of the type $a + b\sqrt{-1}$ various algebraic and transcendental operations which were known in his time, and showed that the result would be a number of the same kind.

### 1752. Jean le Rond d'Alembert. Cauchy–Riemann equations

In the 18[th] century, hydrodynamics developed swiftly. In 1752, d'Alembert examined the perfect fluid motion. In the article *Essai d'une nouvelle théorie de la résistance des fluides*, d'Alembert determined the speed $f(x,y) = u(x,y) + v(x,y)\sqrt{-1}$, where the functions $u(x,y)$ and

$v(x, y)$ are projections of the speed of a particle of liquid on the axis of coordinates. They are connected by the equations $\dfrac{\partial u}{\partial x} = \dfrac{\partial v}{\partial y}$, $\dfrac{\partial v}{\partial x} = -\dfrac{\partial u}{\partial y}$, i.e. $vdx + udy$ and $udv - vdy$ are full differentials, with the compact notation $\dfrac{\partial f}{\partial x} + i\dfrac{\partial f}{\partial y} = 0$. In 1755, Euler arrived at the same results, and later established that the actual and imaginary part of any analytical function necessarily satisfy these conditions:

Figure 5, Euler L. Ulterior disquisitio de formulis integralibus imaginariis[22].

Euler's works set out the theory of elementary functions of a complex variable. Today, these are called Cauchy–Riemann conditions, and are the conditions for the analysis of a function. For this function the family of curves $u(x, y) = C$ and $v(x, y) = C$ are mutually orthogonal.

### 1768 Euler. *Universal Arithmetic*

Operations of deriving the root continued to present difficulties for a long time. In the *Universal Arithmetic* of 1768 in Russian[23], Euler writes "Roots from negative numbers are not more or less than nothing, and they are also not nothing, for 0 multiplied by 0 in derivation gives 0, and accordingly is not a negative number.

When all possible numbers which can be imagined are more or less than 0 or 0 itself, it can be seen that the square roots of negative numbers cannot be included in the group of possible numbers, and accordingly they are *impossible* numbers. This fact leads us to a recognition of these numbers, which by their property are impossible and are usually called *imaginary* numbers, because they can only be imagined in the mind"[24] (Euler's italics).

But Euler goes on to make a mistaken argument: "But when $\sqrt{a}$, multiplied by $\sqrt{b}$, gives $\sqrt{ab}$; then $\sqrt{-2}$, multiplied by $\sqrt{-3}$, will give $\sqrt{6}$; equally $\sqrt{-1}$, multiplied by $\sqrt{-4}$, will give $\sqrt{4}$, i.e. 2; from this we can see that two impossible numbers multiplied together we can get a possible or real number. But when $\sqrt{-3}$ is multiplied by $\sqrt{+5}$, we get $\sqrt{-15}$, or a possible number, multiplied by an impossible one always give an impossible number"[25]. As we can see, operations on complex numbers were not yet clear, but 9 years later Euler corrected his error, defining $\sqrt{-1}$ as the imaginary $i$, the square of which is equal to $-1$, i.e. $\frac{1}{i} = -i$. Euler discussed the issue on the expedience of imaginary numbers: "It finally remains to dispel the doubt about when these forces are impossible, it seems that they are completely unnecessary, and this science may be considered worthless. But despite this, it is in fact very necessary, for such issues often arise, in which we cannot discover swiftly whether they are possible or impossible. But when their solution brings us to these impossible numbers, this will mean that the actual issue is impossible. To clarify this with an example, let us examine the following issue: the number 12 divided into two parts, the product of which would be 40. When we solve this issue in the next rules, we will find for the two numbers $6 + \sqrt{-4}$ and $6 - \sqrt{-4}$, which accordingly are impossible: thus, from this we see that this problem cannot be solved.. If the number 12 must be divided into two such parts, which would give 35, then these parts would undoubtedly be 7 and 5"[26].

### 1777, Euler introduces the symbol *i*

In the paper *On forms of differentials of angles, especially with irrationalities, which are integrated with the assistance of logarithms and circular arcs*, *Master of natural sciences of the Academy presented on 5 May 1777,* published in 1794, for the first time Euler introduced the symbol of the imaginary unit *i*, from the first letter of imaginaire, which is what Descartes called imaginary numbers.

Figure 6. The first appearance of symbol "*i*"

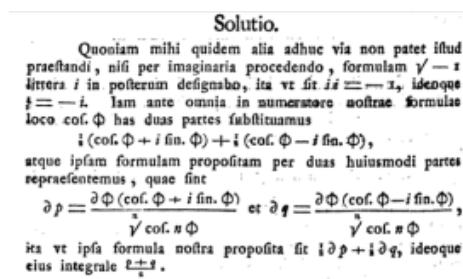

*Translation*: We will examine and study the differential formula $\frac{\partial\Phi\cos\Phi}{\sqrt[n]{\cos .n\Phi}}$, the integral of the logarithm of the circular arcs. Solution. For this I believe another method is also available, which however requires the imaginary unit $\sqrt{-1}$, which in future we will designate by the letter *i, so* $ii = -1$, or $\frac{1}{i} = -i$ .. Above all, we note that the value of our formula, $\cos.\Phi$ can be replaced by two parts $\partial p = \frac{\partial \Phi(\cos.\Phi + i\sin.\Phi)}{\sqrt[n]{\cos.n\Phi}}$ and $\partial q = \frac{\partial \Phi(\cos.\Phi - i\sin.\Phi)}{\sqrt[n]{\cos.n\Phi}}$, and then our formula may be represented as $\frac{1}{2}\partial p + \frac{1}{2}\partial q$, and the integral is expressed as $\frac{p+q}{2}$ [27].

[25] Ibid, p. 93.
[26] Ibid, p. 94–95.

In the 1770s, Euler's changed his attitude towards imaginary numbers. From auxiliary formalism, it acquired the necessary theoretical status, receiving the definition ($i^2 = -1$) and a description of properties. Euler developed the theory of integrals of the function of the complex variable, and also singled out the *principle of symmetry*. "The entire theory of imaginary numbers, to which analysis is now obliged for so much success, rests primarily on the following foundation: if $Z$ is any function from $z$, which after the substitution $z = x + y\sqrt{-1}$ takes the following form: $M + N\sqrt{-1}$, which by the substitution $z = x - y\sqrt{-1}$ the same function $M - N\sqrt{-1}$, where the letters $M$ and $N$ always mean real quantitates"[28]. From this, the Euler–d'Alembert formulas follow, or as we call them today, the Cauchy–Riemann formulas.

Figure 7. Euler on the symmetry property

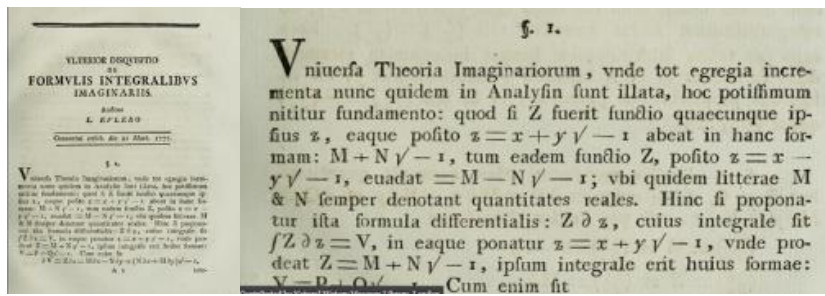

Euler then applied the function of the complex variable to conformal transformations (preserving the angles and likeness in the small).

## 1797/1799, Caspar Wessel

The geometric interpretation of complex numbers and operations with them was first given by the Norwegian geodesist and cartographer of the Danish academy of sciences Caspar Wessel (1745–1818) in the work *An essay on the analytical representation of direction and its*

*applications, primarily to the solution of flat and spherical polygons* [29], submitted in 1797 and published in Danish in 1799. He wrote the work for cartographers.

Wessel introduced the concept of the directed section, and defined structure as the parallel displacement of a plane, and multiplication as the rotation of a plane with expansion. In of his work, Wessel writes: "Let +1 indicate a positive linear unit, and $+\varepsilon$ another different unit perpendicular to the positive unit and having the same origin. Then the directing angle +1 will be equal to $0°$, for −1 will be equal to $180°$, for $+\varepsilon$ will be equal to $90°$, for $-\varepsilon$ will be equal to $-90°$ or $270°$. Owing to the rule that the directing angle of product is equal to the sum of the angle of co-multipliers, we will obtain: $(+1)(+1)=+1$, $(+1)(-1)=-1$, $(-1)(-1)=+1$, $(+1)(+\varepsilon)=+\varepsilon$, $(+1)(-\varepsilon)=-\varepsilon$, $(-1)(+\varepsilon)=-\varepsilon$, $(-1)(-\varepsilon)=+\varepsilon$, $(+\varepsilon)(+\varepsilon)=-1$, $(+\varepsilon)(-\varepsilon)=+1$, $(-\varepsilon)(-\varepsilon)=-1$. From this it is clear that $\varepsilon$ is equivalent to $\sqrt{-1}$ and the deviation of the product is determined in such a way that not one of the general rules violates this operation"[30]. Wessel showed that complex numbers represented by directed segments obey non-contradictory arithmetic. The sum of the two complex numbers $a+bi$ и $c+di$ Wessel calls the diagonal of a parallelogram, built on the sides of directed segments, corresponding to the components, i.e. parallel displacement of the plane along $a+bi$. Multiplication of the two complex numbers $(a+bi)(c+di)=(a+bi)\rho e^{i\varphi}$, where $\rho e^{i\varphi}=c+di$ reflects the rotation of plane around the point $O$ to the angle $\varphi$ with the extension of all scales in the relation of $1:\rho$. Wessel's work contained the foundations of vectoral calculation for two-dimensional space, and was the geometric model of complex numbers, but unfortunately it passed unnoticed in both Denmark and the rest of Europe. Europeans did not read it, as they did not know Danish, and Danish academicians ignored it. Only a century later, in 1897, in Copenhagen a French translation of the work was published, edited by Zeiten. Now it is available in English in Smith's anthology[31]. Wessel's discovery had no influence on European mathematics. In the 19[th] century, the geometric interpretation of the complex number was once more discovered by Argan, and developed in the works of Gauss, Grassmann, Hamilton and other scientists.

### 1806, 1813/14. Jean-Robert Argand (1768–1822).

In 1806, the bookseller[32] Jean Robert Argand (1768–1822) anonymously published a brochure *An essay on a certain method of representing imaginary values in geometric structures*

---

[33]. Argand developed a geometric theory of the complex number, drawing the same conclusions as Wessel. In particular, he noticed that in the multiplication of complex numbers their arguments combine (Argand, p. 20), and the modules stretch. Argand introduced the so-called Argand diagrams, depicting multiplication operations, raising to a power and extracting a root from a complex number.

Argand came close to the concept of trigonometric polynomials of Chebyshev (least deviating from zero) $\cos 2\alpha = 2\cos^2\alpha - 1$, $\sin 2\alpha = 2\sin\alpha\cos\alpha$, $\cos 3\alpha = 4\cos^3\alpha - 3\cos\alpha$, $\sin 3\alpha = 3\sin\alpha - 4\sin^3\alpha$ etc., i.e. $F_1(x) = x$, $F_2(x) = 2x^2 - 1$, $F_3(x) = 4x^3 - 3x$ etc.

Argand's works was republished in 1813/14 by G. Gergonne in the 4[th] volume of the journal *Annales de mathématiques pures et appliquées* (volume 4, 1813–1814, p. 61) together with his new article (ibid, p. 133) in the 5[th] volume (Volume 5, p. 197). Other works also appeared on this topic, and Hankel writes about them.

Figure 8.Here is an example of one of Argand's diagrams from the English edition[34] :

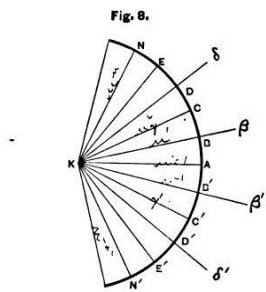

Fig. 8.

§1. If AB, BC, . . . . , EN (Fig. 8) are equal arcs, $n$ in number, and we make $\overline{KB} = u$, we shall have $\overline{KC} = u^2$, $\overline{KD} = u^3$, . . . . . , $\overline{KN} = u^n$.

## 1821 A. Cauchy (1789–1857). *Analyse algébrique*

In 1821, Augustin Cauchy delivered a course of analysis at the Polytechnical School. As is well known, his students protested vehemently against the instruction of complex numbers, believing this knowledge to be useless. The description of this topic in *Analyse algébrique* is formal: Cauchy examines operations on complex numbers as operations on algebraic symbols, which was later mocked by Hankel: addition theorems, trigonometric form, arithmetic operations on algebraic form, reduction of algebraic form to trigonometric and vice versa; raising to a power and extraction of a root. However, Cauchy was the first to give the formula $\left(\cos\vartheta + \sqrt{-1}\sin\vartheta\right)\left(\cos\vartheta - \sqrt{-1}\sin\vartheta\right) = 1$, and also to notice the periodicity of the complex

---

number[35]. Cauchy does not provide a geometric interpretation of complex numbers and operations with them. But this was a curriculum, not a scientific study. Later, in 1829–32 Cauchy made a great contribution to the theory of functions of the complex variable – he created the Residue Theory.

### 1831 Carl Friedrich Gauss. *Theory of Biquadratic Residues*

Gauss (1777–1855) was one of the most secretive scientists of his time. He knew everything, but published very little. An understanding of the geometric nature of complex numbers is latently contained in his dissertation of 1799, but a strict establishment of the algebra of complex numbers was made in *Theory of biquadratic residues* of 1831. Gauss wrote: "The difficulties that are thought to surround the theory of imaginary values largely arise from their rather unfortunate names (some have even given them the unfortunate sounding name of impossible values). If judging from the concepts given by the diversity of two dimensions (which with great clarity are manifested in spatial considerations) to call positive values direct, negative reverse, and imaginary perpendicular values, we would have simplicity instead of confusion, clarity instead of vagueness"[36]. He examined numbers on a complex plane, introduced the concept of the adjoined number, the norm, and classified whole and rational complex numbers. As Gauss' goal in this work was number theory, he also brought in complex numbers with this goal. As some simple numbers proved to be factors of conjugate complex numbers, for example, $2 = (1+i)(1-i)$, $5 = (1+2i)(1-2i)$, $13 = (3+2i)(3-2i)$, $17 = (1+4i)(1-4i)$, Gauss was able to discover the theoretical-numerical patterns of the factorization.

### 1841, H. Grassmann. *The Theory of Linear Extension*

Hermann Grassmann (1809–1877), a teacher at a Prussian grammar school, studied the nature, configuration and multiplication of complex numbers, and on this basis wrote the *The Theory of Linear Extension, a New Branch of Mathematics* (*Ausdehnungslehre*), published in 1844 and subsequently in 1862. Based on the principles and "requirements of statics and mechanics"[37], he introduced the concept of *n*-dimensional manifold with a system of operations. In particular, Grassmann, wrote: *"By the product of two segments a,b, we understand the area of the parallelogram that they form, meaning both the value and its position*, i.e. we assume $ab = cd$ only in the case if the parallelogram formed by the segments *a* and *b* is not only equal in size to the parallelogram that is formed by the segments $c$ and $d$, but which also lies on a plane parallel with the latter, and has the same direction[38] (Italics in original) <...> If we change the place of the factors of product of $ab$, then the sense of the parallelogram changes to the reverse" [39]. But the difficulty of exposition and the philosophical language in place of mathematical language made it difficult to understand his discovery for a long time, which only by the end of the 19[th] century served as the foundation for introducing the $n$-dimensional vectoral space made in the works of Gibbs.

### 1843, W.R. Hamilton. Creation of the theory of Quaternions

Sir William Rowan Hamilton (1805–1865), was the royal astronomer of Ireland, a mathematician, theoretical mechanic and theoretical physicist. From 1835, Hamilton began to regard algebra not as an art, language, or science of quantity, but rather as a science of order in certain rows. An example of this process was for him ideal time, freed of all connections of cause and effect, as according to Kant it is the pure intuitive form of our inner perception, and therefore better adapted than space, i.e. the form of our external perception; at any rate, the concept of "past", "present" and "future" arise in our consciousness sooner than the concepts of "forward" and "backward" in space; therefore, algebra for him is a science of pure time. "If geometry relies on intuition of space, algebra could rely on related intuition of time"[40]. Hamilton determined vector as shift. His symbol $i$ means firstly the single vector of the axis $Ox$, secondly an imaginary unit, and thirdly the operator of rotation – the versor.

In 1835, Hamilton published the work *The Theory of Algebraic Pairs* [41], in which he gave a new structure of the theory of complex numbers. This was the following form $\left(z=\left(a,b\right)\right)$ of complex numbers after the algebraic $\left(z=a+bi\right)$, trigonometric $\left(z=r\left(\cos\varphi+i\sin\varphi\right)\right)$ and exponential $\left(z=re^{i\varphi}\right)$. Hamilton examined the complex number $x+iy$ as an algebraic pair $\left(x,y\right)$ of real numbers, i.e. the removed the geometric element and reduced complex numbers to pure algebra, making it possible to move to a new level of geometric generalization – turn and extension on a plane. This formalized methods of mathematical physics in tasks of the flow of liquid or heat, gravitation, sound and optics. But these tasks were solved in a two-dimensional space.

---

Hamilton wished to apply the system of complex numbers to three-dimensional space, but found difficulties with determining multiplication – either the commutative law or the law of distribution was violated. This contradicted the principle of permanency of the equivalent forms of J. Peacock, established in 1830: laws of operations of algebra must remained unchanged, whatever the symbols mean on which operations are carried out[42]. Furthermore, for Hamilton non-zero multipliers could give a zero product. Hamilton conclude that algebra could only be organized for four-dimensional numbers.

Hamilton had a sudden inspiration of how to multiply four-dimensional numbers on the 16th of October 1843 in Dublin, on Broom Bridge across the Royal canal: «And here there dawned on me the notion that we must admit, in some sense, a fourth dimension of space for the purpose of calculating with triplets; or transferring the paradox to algebra, must admit a third distinct imaginary symbol $k$, not to be confounded with either $i$ or $j$, but equal to the product of the first as multiplier, and the second as multiplicand; and therefore was led to introduce quaternions, such as $a+ib+jc+kd$, or $(a,b,c,d)$.»[43]. Hamilton was so amazed that he immediately scratched the formulas on the rails: $i^2 = j^2 = k^2 = ijk = -1$.

Figure 9. A memorial plaque on Broom Bridge commemorates his discovery.

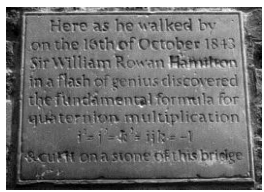

As Hamilton himself wrote in 1843 about rotation "My train of thoughts was of this kind. Since $\sqrt{-1}$ is in a certain well-known sense, a line perpendicular to the line 1, it seemed natural that there should be some other imaginary to express a line perpendicular to the former; and because the rotation from this to this also being doubled conducts to −1, it ought also to be a square root of negative unity, though not to be confounded with the former. Calling the old root, as the Germans often do, $i$, and the new one $j$, I inquired what laws ought to be assumed for multiplying together $a+ib+jc$ and $x+iy+jz$.»[44].

For the "four dimensional numbers" that he discovered, Hamilton introduced the name quaternions – from Latin quaterni, "four of each". He wrote out quaternions as the sum of the type $q = a + bi + cj + dk$, where $\boldsymbol{i, j, k}$ are three quaternion units (equivalents of the imaginary unit *i*), and *a*, *b*, *c*, *d* are real numbers. By assuming the multiplication of quaternions to be distributive in relation to structure, Hamilton reduced the definition of the operation of multiplying quaternions to the task of the table of multiplication for the basic units 1, $\boldsymbol{i}, \boldsymbol{j}, \boldsymbol{k}$ :

| × | 1 | $i$ | $j$ | $k$ |
|---|---|-----|-----|-----|
| 1 | 1 | $i$ | $j$ | $k$ |
| $i$ | $i$ | $-1$ | $k$ | $-j$ |
| $j$ | $j$ | $-k$ | $-1$ | $i$ |
| $k$ | $k$ | $j$ | $-i$ | $-1$ |

From the table we can see that the multiplication of quaternions is not commutative (so the algebraic system of quaternions is not a field). The addition of vectors is commutative, it means parallel space transfer. But the result of implementing two three-dimensional turns in multiplication depends on their order. Rotation around the origin in space defines an axis, and so the extension with the rotation, which in the case of the plane required two constants, in space can be characterized by only four parameters.

The theory of vectoral function of the scalar argument was developed by Hamilton in 1846. Hamilton also developed the concepts of collinearity and coplanarity of vectors, orientation of vectoral trio and others. He introduced the concept of the hodograph, the nabla operator, and applied the theory to tasks of celestial mechanics.

### The appearance of vectoral analysis

Hamilton's theory of quaternions, which he described in 109 articles, was summarized concisely by his pupil, the Scottish mathematician Peter Guthrie Tait. A friend and fellow pupil of Tait, J.C. Maxwell (1831–1879), saw in the theory of quaternions a convenient apparatus for the mathematical description of the theory of electricity and magnetism, contained in the concept of field and power lines, described by Michael Faraday in 1839–1855. Maxwell separated vectoral calculus from the theory of quaternions. He did so in the work *A Treatise on Electricity and Magnetism* (1873)[45] in the section "Preliminary". Maxwell's work contains almost no symbols of quaternions, but the most useful for tasks of physics is taken from it. Maxwell called the vector $-\left( i\dfrac{d\psi}{dx} + j\dfrac{d\psi}{dy} + k\dfrac{d\psi}{dz} \right)$ a *slope* of the function ψ, to show the direction and rate of fastest decrease of ψ [46], and for the function of two variables, the direction of the steepest slope

---

of the surface. The term *gradient* is derived from the Latin *gradior* – "move ahead". The term became used in meteorology, and Maxwell later used it to replace his slope of $\psi$.

Over time, the square of the imaginary value $i^2 = -1$ was replaced by the scalar product $(i,i) = 1$. Later, Gibbs wrote *Elements of vector analysis*[47] (1880s), after which Heaviside (1903) gave vectoral calculus its modern form. The term "vector analysis" was proposed by Gibbs (1879) in his course of lectures. Gibbs' description of vector analysis became classic[48].

Quaternions are still used in geometry and physics, for example in the Lorentz transformation, where it is important to set a three-dimensional turn with the help of a minimum number of scalar parameters, this description never degenerates.

In a report made on 1 September 1908 at the 80[th] meeting of German scientists and doctors in Cologne, H. Minkowski called the totality of the substantial *x, y, z, t* a *world* [49].

### 1867. Hermann Hankel. *Theory of complex number systems*

In 1867, Hermann Hankel (1839–1873) published his summarizing book *Theory of complex number systems, primarily ordinary imaginary numbers and Hamilton quaternions together with their geometric interpretation by Doctor Hermann Hankel* (*Vorlesungen über die komplexen Zahlen und ihre Functionen*, 1. Teil. *Theorie des complexen Zahlensysteme*).

Hankel gives a historical survey and analysis of complex numbers and their systems. For example, where the law of commutativity is observed, but associativity and distributivity is violated (Scheffler's work of 1851[50]), Kirkman's systems of complex numbers[51], does not obey the law of associativity, connected with the operational hyperdeterminant of Arthur Cayley. Hankel notes that "there is a connection between function theory and complex numbers of the higher order and so-called operational calculation, simply the symbolic combination of certain operations on numbers"[52]. It was thanks to Hankel's suggestion that Grassmann's work *The Theory of Linear Extension* became understood and gained recognition.

Hankel advanced the hypothesis: "A chemical formula may be examined like a complex number, the units of which are served by the chemical designations of elements, and the coefficients signs showing the total number of each of the elements. A chemical compound corresponds in number theory the operation of multiplication; elements or their atomic weights correspond to initial multipliers, and chemical formulas for the decomposition of bodies are literally the same as formulas for the decomposition of numbers"[53]. This was two years before Mendeleev discovered the periodic table.

Hankel formulated the law of permanence of formal laws: "If two parts of a logical form expressed by common signs of universal arithmetic are equal, then they should also remain equal when the signs that express them cease to indicate ordinary values, and as a consequence the operations themselves gain a somewhat different, but definite meaning" [54].

## Conclusion

Real numbers, the theory of which differed from the method of exhaustion of Eudoxus before the concept of the uninterrupted numerical region created simultaneously by Méray, Heine, Cantor, Dedekind and Weyerstrass, was summarized by A.N. Kolmogorov[55]. Real numbers are closed relative to arithmetic operations, and regulated.

The appearance of an imaginary unit broadened the multitude of real numbers, forming a two-dimensional space – a complex plane. This is also a complete and single system, but it no longer has any order in it.

Quaternions have a dimensionality of 4 and lose commutativity of multiplication. In 1898, Adolf Hurwitz proved that the system of quaternions is also single. Quaternions are the only finite algebra with division, which contains real numbers and does not coincide with real or complex numbers[56].

The concept of the complex number was developed based on the inner logic of mathematics, and also based on the requirements of applied sciences – cartography, hydrodynamics and other natural sciences. Studies showed that operations on complex numbers reflect the properties of movement in space – rotation and extension.

Gradually, the picture of the world changed: Newtonian mechanics were replaced by relative, and geometric concepts of space also changed. The mathematic axioms of metric space were formed in the first decades of the 20[th] century (Fréchet and Hausdorff). The imaginary component of complex functions gained the physical meanings of projection of force, reactive resistance, energy loss, comprising the refraction ratio, and harmonic oscillation.

---

Here we have broadly discussed the history of how the concept of complex numbers developed. There were also other attempts to interpret complex numbers and operations on them, which did not gain recognition. In the late 19[th] century, complex multi-dimensional systems arose in attempts to expand the concept of complex number, but discussion of this topic goes outside the boundaries of this article.

The development of the concept of complex numbers was the basis for the theory of functions of the complex variable, operational methods of solving differential equations, non-Euclidian geometry, number theory, geodesy and cartography, mathematical physics, elasticity theory, electromagnetic field theory, electro and magneto static, electro dynamics and quantum mechanics. The complexity reflects the fundamental qualities of the world – symmetry and cyclicity.